\documentclass{amsart}
\usepackage[T1]{fontenc}
\usepackage{amssymb}
\usepackage{graphicx}
\usepackage{comment}
\newtheorem{theorem}[subsection]{Theorem}
\theoremstyle{definition}
\newtheorem{definition}{Definition}[section]
\theoremstyle{lemma}
\newtheorem{lemma}{Lemma}[section]
\theoremstyle{program design}
\newtheorem{program design}{Program Design}[section]

\title{Nondeterministic Infinite Time Turing Machines}
\author{Erin Carmody}
\date{}

\begin{document}
\maketitle

\begin{abstract}
This paper analyzes infinitary nondeterministic computability theory.  The main result is D $\ne$ ND $\cap$ coND where D is the class of sets decidable by infinite time Turing machines and ND is the class of sets recognizable by a nondeterministic infinite time Turing machine.  Nondeterministic infinite time Turing machines are introduced, along with relevant definitions.  
\end{abstract}

\section{Introduction}

This paper introduces a concept of nondeterministic computation with infinite time Turing machines (ITTMs) and analyzes the resulting nondeterministic computability theory.  Nondeterministic infinite time Turing machines (NITTMs) are programs which are not well-defined functions, as are deterministic programs, from a given state and cell value to another state, cell value, and instruction for how to move the head.  NITTM programs may respond to a current state and cell value with any number of instructions to be performed simultaneously.  The mechanisms of the machine are the same as an ITTM as detailed in [6].  Let $A = \{ (s,b) \ : \ s$ is a valid state, $b$ is a cell value$\}$ and $B = \{(s,b,i): s$ is a valid state, $b$ is a cell value, $i$ is a head movement instruction$\}$.  Then an ITTM program is a function $f:A \to B$, while an NITTM program is strictly a relation $R$ on $A \times B$.  \\

To better understand how an NITTM program works and to demonstrate the power of a simple NITTM program, allow me to describe a particularly helpful program.  Beginning in the starting state on a machine with input $0$, $(s,0)$, \\

\begin{center}

  \includegraphics[width=.5 \textwidth]{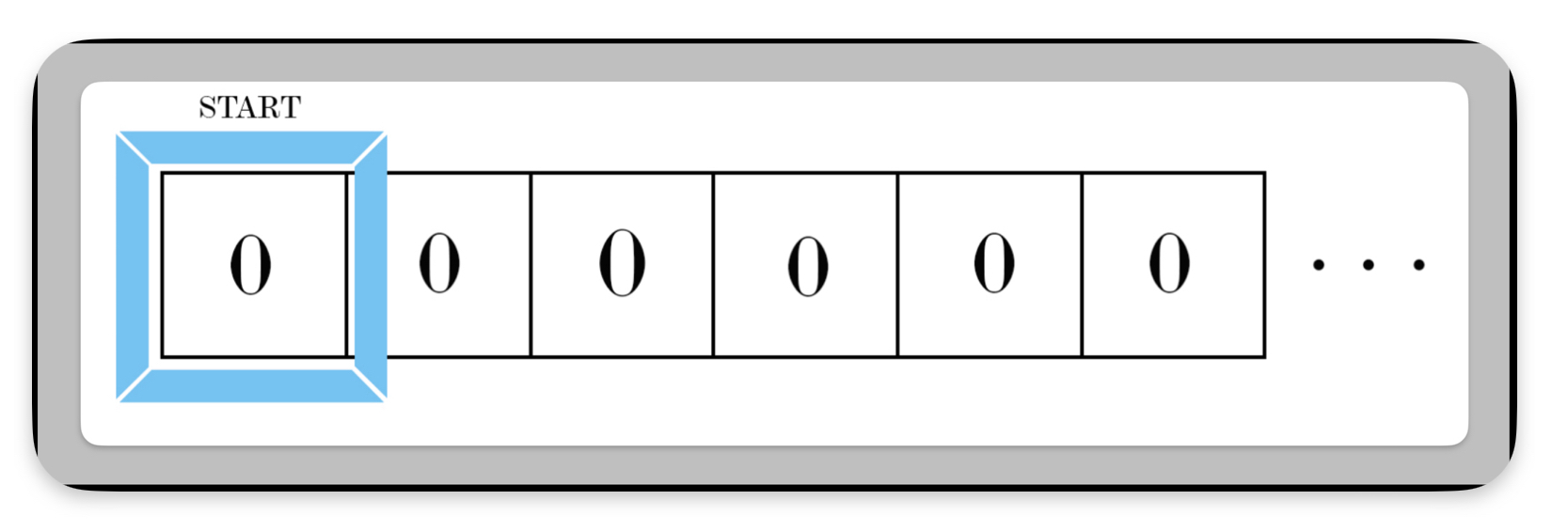}
\end{center}

the program branches with two diverging instructions.  The first is to stay in the start state, write a $0$ on the cell, and move the head to the right $(s, 0, R)$.  The second is to stay in the start state, write a 1 on the cell, and move the head to the right, $(s,1,R)$.  In both cases, the program is again in the start state on the second cell reading a 0, thus in $(s,0)$, and so the program responds again with both $(s,0,R)$ and $(s,1,R)$.  \\

\begin{center}

  \includegraphics[width=.8 \textwidth]{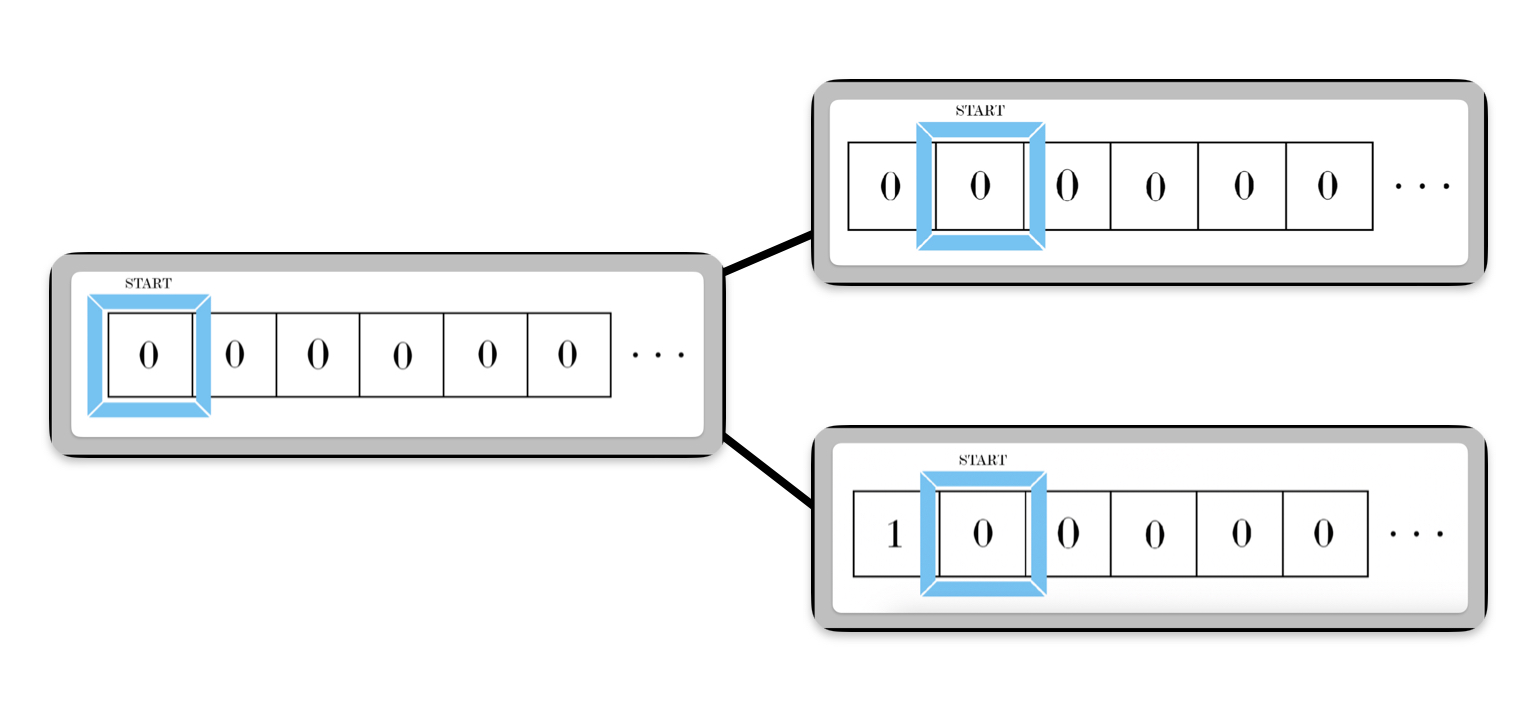}
\end{center}

This will continue for $\omega$ many steps until the machine goes into the limit state $L$, and the head is back at the left-most cell. The cell value of the first cell will be, simultaneously, 0 and 1. 

\begin{center}

  \includegraphics[width=.29 \textwidth]{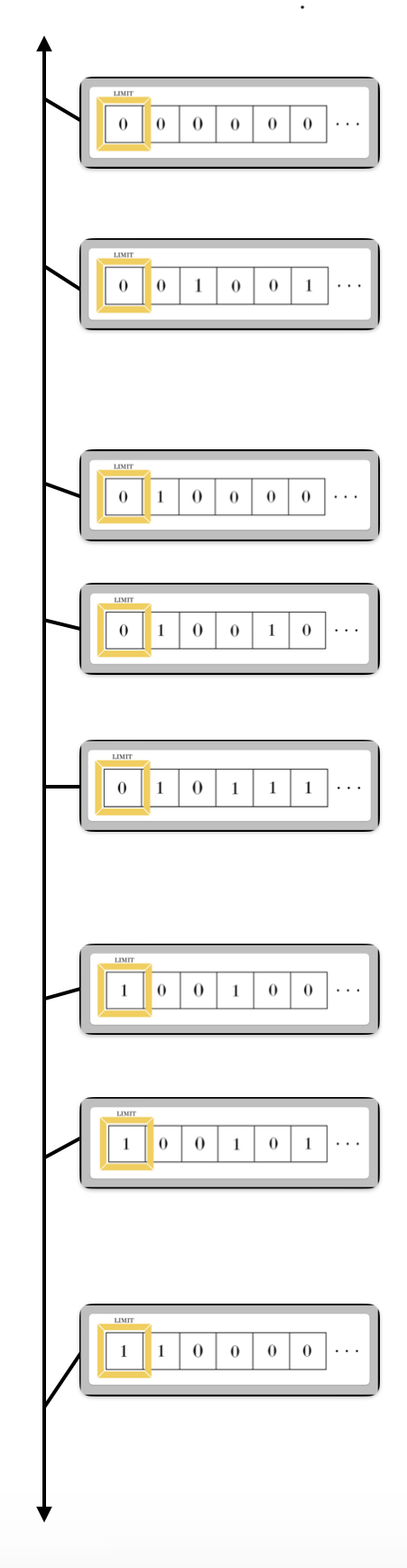}
\end{center}

 Let the program respond to either $(L,0)$ or $(L,1)$ by halting. 
 
 \begin{center}

  \includegraphics[width=.7 \textwidth]{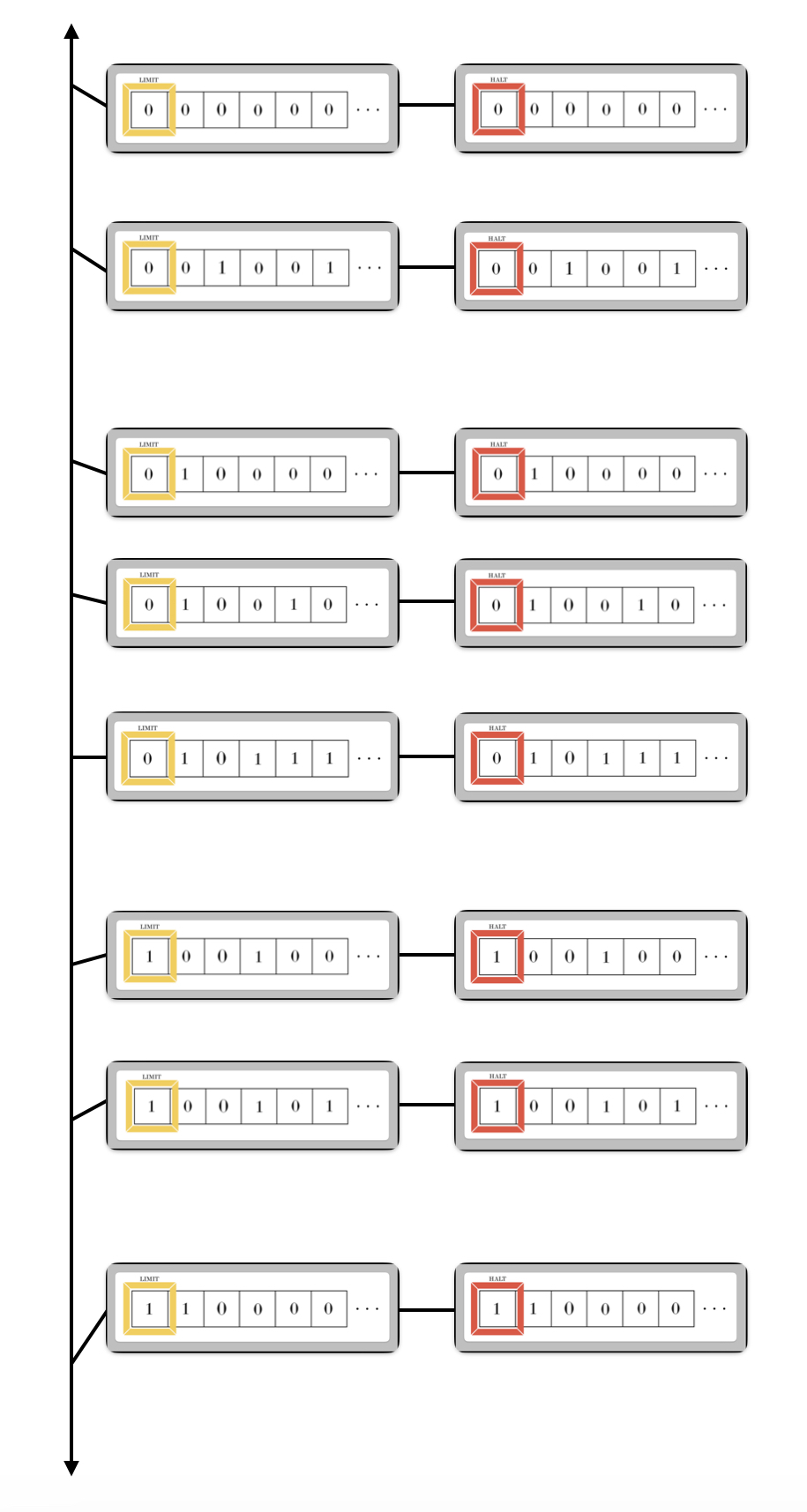}
\end{center}

 The program has nondeterministically written all the real numbers on the tape in $\omega$ many steps.  Thus, the computation tree of the program includes the full binary tree. For future reference, call this program $p_{\mathfrak{c}}$, which will also refer to such a program with instructions to write the values on the scratch tape or output tape instead.\\

\begin{definition}
A program \textit{branches} when for a single given input (state, cell value), there are multiple different program instructions for output (state, cell value, head movement).  \\
\end{definition}

\begin{definition}
A \textit{snapshot} is a list of data $(p,t,k,s)$ where $p$ is a program, $t$ is an infinitary binary sequence describing the contents of the tape, $k$ is a natural number giving the position of the head, and $s$ is the current state of program $p$.\\
\end{definition}

\begin{definition}
A \textit{computational sequence} is a list $\langle x_{\alpha} \ : \ \alpha < \beta \rangle$, for some ordinal $\beta$, such that $x_0$ is the snapshot of the input configuration, each $x_{\alpha + 1}$ is a snapshot obtained from $x_{\alpha}$ as a program response to $(s_{\alpha}, t_{\alpha}(k_{\alpha}))$.  For limit ordinals $\lambda$, the snapshot $x_{\lambda}$ is obtained from the snapshots $\langle x_{\alpha} \ : \ \alpha < \lambda \rangle$ by adhering to limit operations of ITTM machines so that $s_{\lambda} =$ limit state, $t_{\lambda} = \limsup_{\alpha < \lambda} t_{\alpha}$, and $k_{\lambda} = 1$.  \\
\end{definition}

\begin{definition}
A \textit{halting} computation sequence is one with a final snapshot at some ordinal $\beta$ with $s_{\beta} =$ halt state.\\
\end{definition}

\begin{definition}
The \textit{computation tree} of a nondeterministic program $p$ on input $a$ consists of all computation sequences for $p$ with input $a$.  If $q$ and $r$ are computation sequences in $p$, then $q < r$ in the computation tree when $q$ is an initial segment of $r$.  A path through the tree then corresponds with an allowed computation sequence of the program $p$ on input $a$.  The diagram below shows the first three stages of the computation tree and the last stages.\\

\begin{center}

\includegraphics[width=1.0\textwidth]{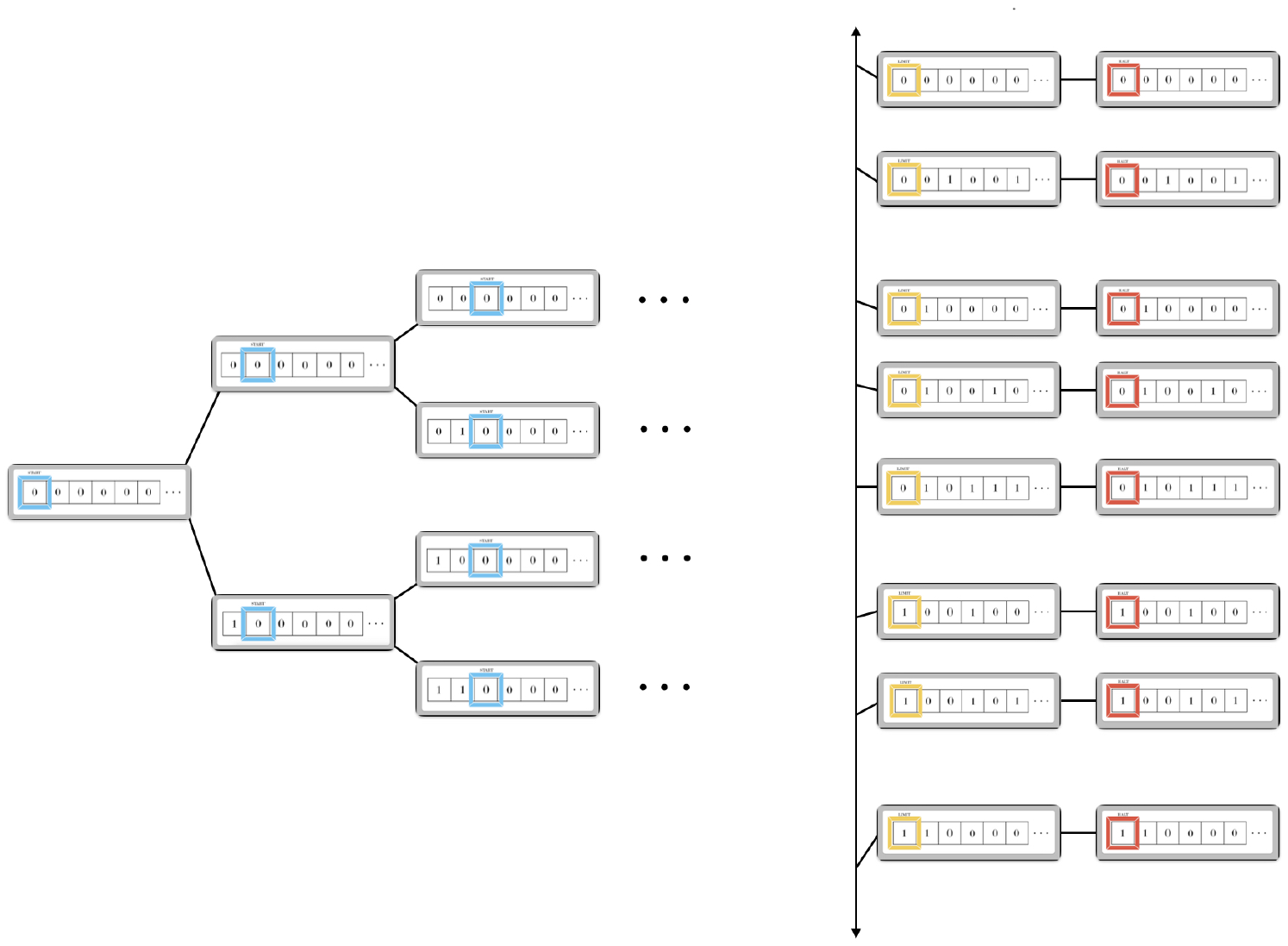}
\end{center}
\end{definition}

Each computation sequence of the computation tree of program $p_{\mathfrak{c}}$ described earlier is of the form $\langle x_n \ : \ n \le \omega \rangle$, where each $t_n \in x_n$ is a finite binary sequence, and $\bigcup_{n \in \omega} t_n$ is an infinitary binary sequence.  So, actually the full binary tree is just a part of the computation tree of $p_{\mathfrak{c}}$.  If $\langle p_{\alpha} \ : \ \alpha < \mathfrak{c} \rangle$ is the collection of computation sequences of $p_{\mathfrak{c}}$, each $p_{\alpha} = \langle x_{n,\alpha} \ : \ n < \omega \rangle$ (the computation sequence minus the halt state snapshot), where $x_{n,\alpha} = (p_{\mathfrak{c}},t_{n,\alpha},k_{n,\alpha},s_{n,\alpha})$ is a snapshot, then the full binary tree is $\bigcup_{n \in \omega} \bigcup_{\alpha < \mathfrak{c}} t_{n,\alpha}$.\\

\section{NITTM Computability Theory}

The main goal of the paper is to analyze NITTMs infinitary computability theory.  In doing so, we meet the P vs. NP problem in an infinite context.  Schindler showed that P $\ne$ NP for infinite time Turing machines [9].  And, it has been shown Deolalikar, Hamkins and Schindler in [3] that the decidable ITTM sets are a proper subclass of $\Delta_2^1$.  Nondeterministic infinite time Turing machines, we shall see, can compute a bit higher. Here I will define the analogous notions of semi-decidability and co-semi-decidability.\\

\begin{definition}
A set $A$ is \textit{recognizable} if and only if there is an NITTM program $p$ such that for all $x \in A$ there is a computation sequence for program $p$ on input $x$ such that the final snapshot of the sequence includes the halt state with output 1.  
\end{definition}

\begin{definition}
The class of sets \textit{ND} are those which are recognizable.\\
\end{definition}

\begin{definition}
A set $B$ is \textit{co-recognizable} if and only if the compliment of $B$ is recognizable.\\
\end{definition}

\begin{definition}
The class of sets $coND$ are those which are co-recognizable.\\
\end{definition}

\begin{definition}
The class of sets $D$ are those which are ITTM decidable.\\
\end{definition}

The following Lemma establishes which sets are recognizable by an NITTM program.\\

\begin{lemma}
A set $A$ is recognizable if and only if $A$ is a set of $\Sigma_2^1$ set of reals.

\begin{proof}
A $\Sigma_2^1$ set of reals is one such that $x \in A \iff \exists s \ \forall t \ \exists u\  B(x,s,t,u)$ where $B$ is a $\Delta_0^1$ set.  Suppose first that $A$ is $\Sigma_2^1$ set.  Then $A$ is as just described, or we can say that  $A$ is of the form $\exists s \ B(t,s)$ where $B$ is a $\Pi_1^1$ set of reals.  It is known by [3] that ITTM programs can decide all $\Pi_1^1$ sets.  Let $p$ be the ITTM program which decides $B$.  For the lemma, we only need $p$ to be the program which semi-decides $B$.  Then for any input $(x,y)$, the program $p$ halts with output 1 if $(x,y) \in B$.  Then, using $p$, design an NITTM program which can recognize $A$ as follows.  The $\omega$th level of the computation tree of this program is the same as the computation tree of the program $p_{\mathfrak{c}}$ with each path in the tree containing a real number.  It is possible that any real number is written on the scratch tape by the $\omega$th step of the program.\\

Next, the program instructions at $(L,0)$ and $(L,1)$ is $(s_1,b_1,i_1)$ where this triple is the output of the program $p$ upon input configuration.  The input of the program is a real number $x$, which is written on the input tape.  There are continuum many computation sequences and continuum many universes of computation.  Each universe has a different real $y$ composed on the scratch tape.  In each world, program $p$ checks whether $(x,y)$ is in $B$.  If any computation sequence is halting with a final snapshot including $t=1$, then conclude $x \in A$.  Thus, $A$ is recognizable by an NITTM program.\\

For the other direction, we show that if $A$ is recognizable, $A$ can be expressed as a $\Sigma_2^1$ statement.  A set $A$ is recognizable if there is a program $e$ such that for all $x \in A$ there is a computation sequence $s$ for program $e$ on input $x$ such that the final snapshot is halting with output 1.  To solve the problem, we need to code $e$ and computation sequence with a real.  An NITTM program $e$ is a finite list of instructions, with the option of only finitely many states.  Thus there are only countably many NITTM programs.  Therefore, we can easily assign a natural number to each program $e$.\\

To code a computation sequence by a real, first start by considering the length of a computation sequence.  Since every ITTM program halts or repeats by stage $\omega_1$, it follows that every halting computation sequence, $s$, is of length $\beta$ for some countable ordinal $\beta$.  Hence the length of any stage of a computation sequence can be coded by a real.  Each term $x_{\alpha}$ in the computation sequence is a 4-tuple of natural numbers an thus can be G\"{o}del coded by a natural number $n$.  Thus, the sequence is coded by $[s] = \langle (n_{\alpha},r_{\alpha}) \ : \ n \in \mathbb{N}, \ \alpha < \beta \rangle$ where $n_{\alpha}$ codes the snapshot at stage $\alpha$ , and $r_{\alpha}$ codes $\alpha$.  Next, the pair $(n_{\alpha},r_{\alpha})$ can be coded by a single real, still call it $r_{\alpha}$.  Then the computation sequence is now coded as $[s] = \langle r_{\alpha} \ : \ \alpha < \beta \rangle$.  Finally, this $\beta$ sequence can be coded by a single real.  Hence, there is an $r \in \mathbb{R}$ which codes the given computable sequence $s$.  Hence, the definition of recognizable, $\exists $ an NITTM program  $e$ such that $\forall x \in A \ \exists$ a halting sequence $s$ for program $e$ on input $x$ with output 1, is a $\Sigma_2^1$ statement, so $A$ is $\Sigma_2^1$.
\end{proof}
\end{lemma}

The discussion of NITTMs includes the yet unsolved P vs. NP problem for finite computations.  For infinite time Turing machines we have discussed that the question is settled with inequality as P $\ne$ NP and further P $\ne$ NP $\cap$ coNP in [9] and [3] respectively.  Fagin showed in 1974 [4] that the class NP is exactly those classes of structures axiomatizable in existential second-order logic on finite structures [1], that is the class $\Sigma_1^1$ on finite structures.  Thus co-NP, is characterized by universal second-order logic, that is the class $\Pi_1^1$ on finite structures.  Though, we do not know if $\Sigma_1^1$ is different from $\Pi_1^1$ on finite structures [1], so we cannot conclude that NP $\ne$ co-NP from Fagin's result.  In the infinitary analogue, we have ND is the class $\Sigma_2^1$, and co-ND is the class $\Pi_2^1$, and so ND $\ne$ co-ND.  Thus D $\ne$ ND and further  D $\ne$ ND $\cap$ coND, as shown in the following theorem.

\begin{theorem}
D $\ne$ ND $\cap$ coND.

\begin{proof}
For the deterministic part of the equation, it is established in [3] that the class of sets decidable by an ITTM is a proper subset of $\Delta_2^1$. So $D \subsetneq \Delta_2^1$.  For the nondeterministic side, Lemma 2.1 proves ND $= \Sigma_2^1$.  It immediately follows that coND $= \Pi_2^1$.  Thus ND $\cap$ coND is the full class $\Delta_2^1$.  Thus D $\ne$ ND $\cap$ coND, in particular D $\subset$ ND $\cap$ coND.  

\end{proof}
\end{theorem}

\vspace{.2 in}

\begin{center}

  \includegraphics[width=1 \textwidth]{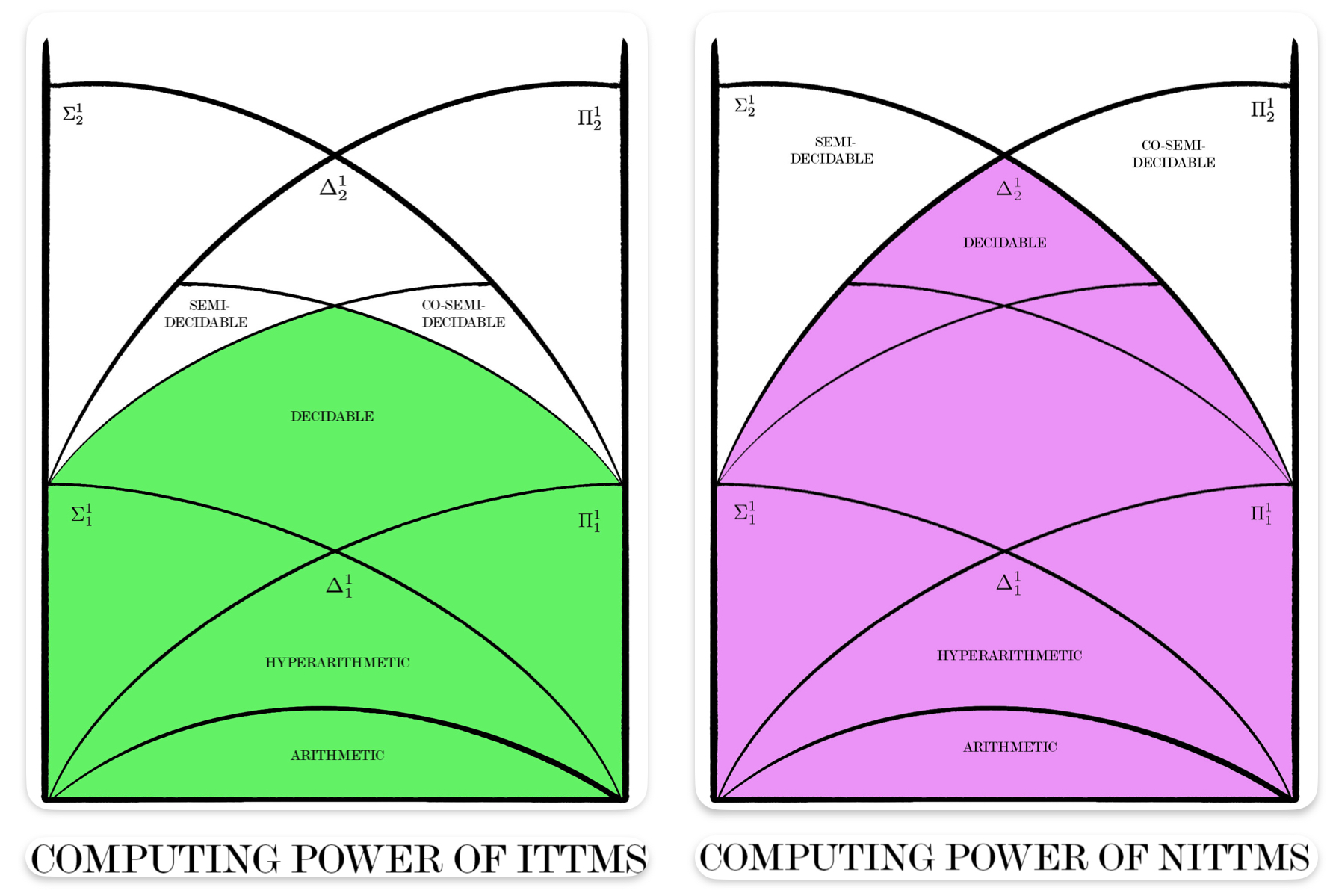}\\

\end{center}
\vspace{.2in}

The same strategy in this proof will not work in the finite case because here we know that ND $\ne$ co-ND because $\Sigma_1^2 \ne \Pi_1^2$, but this is not necessarily true for $\Sigma_1^1$ and $\Pi_1^1$ on finite structures [1].  Also, in the finite case we know that the computing power of a polynomial time Turing machine is a subset of $\Delta_1^1$ on finite structures, but do not know if it is a proper subset [10].  In the infinite case we know that the computing power of an infinite time Turing machine is a proper subset of $\Delta_2^1$.  The figures above are a recreation of the original diagram in [6] with the computing power of ITTMs highlighted, along with the diagram highlighting the analogous computing power of NITTMs. \\

\end{document}